\documentclass[12pt,twoside]{article}

\usepackage{amsmath, amssymb}
\usepackage{graphicx}

\usepackage{color}

\allowdisplaybreaks 

\newtheorem{theorem}{Theorem}[section]

\newtheorem{example}[theorem]{Example}

\newtheorem{lemma}[theorem]{Lemma}

\newtheorem{remark}[theorem]{Remark}

\newtheorem{Thm}[theorem]{Theorem}

\newenvironment{proof}[1][Proof]{\textbf{#1.} }{\ \rule{0.5em}{0.5em}}

\def\vare{{\varepsilon}}

\def\Ad2{{\| A-\widetilde A \|_{L^{2}_{loc}} }}

\def \eref#1{\hbox{(\ref{#1})}}

\def\th{{\theta}}

\title{Least squares estimators for discretely observed stochastic processes driven by small L\'{e}vy noises}
\author{Hongwei Long$^{a,}$\footnote{Corresponding author. Tel: +1 561 2970810; fax: +1 561 2972436. E-mail address: {\tt hlong@fau.edu}}, Yasutaka Shimizu$^{b}$ and Wei Sun$^{c}$ \\
{\it \small $^{a}$Department of Mathematical Sciences, Florida, Atlantic University} \\
{\it \small Boca Raton, Florida 33431-0991, USA } \\
{\it \small $^{b}$Graduate School of Engineering Science, Osaka University} \\
{\it \small Toyonaka, Osaka 560-8531, Japan} \\
{\it \small $^{c}$Department of Mathematics and Statistics, Concordia University}\\
{\it \small Montreal, Quebec H3G 1M8, Canada} \\
}
\date{May 18, 2012}

\begin{document}

\maketitle
\begin{abstract}
We study the problem of parameter estimation for discretely observed stochastic processes driven by additive small L\'{e}vy noises. We do not impose any moment condition on the driving L\'{e}vy process. Under certain regularity conditions on the drift function, we obtain consistency and rate of convergence of the least squares estimator (LSE) of the drift parameter when a small dispersion coefficient $\varepsilon \to 0$ and $n \to \infty$ simultaneously. The asymptotic distribution of the LSE in our general setting is shown to be the convolution of a normal distribution and a distribution related to the jump part of the L\'evy process.

\begin{flushleft}
\vspace{1mm}
\end{flushleft}
\end{abstract}

\section{Introduction}\label{sec:intro}

Let  $(\Omega, {\cal F}, {\mathbb P})$ be a basic probability space
equipped with a right continuous and increasing family of
$\sigma$-algebras $({\cal F}_{t}, t\geq 0)$.
Let $(L_{t}, t\geq 0)$ be a ${\mathbb R}^{d}$-valued L\'{e}vy process, which is given by
\begin{equation}
L_{t}=at+\sigma B_{t}+\int_{0}^{t}\int_{|z|\leq 1} z\tilde{N}(ds,dz)+\int_{0}^{t}\int_{|z|> 1} zN(ds,dz),
\label{e.1.1}
\end{equation}
where $a =(a_{1},\dots,a_{d})\in {\mathbb R}^{d}$, $\sigma=(\sigma_{ij})_{d\times r}$ is a $d\times r$ real-valued matrix, $B_{t}=(B_{t}^{1},\dots,B_{t}^{r})$ is a $r$-dimensional standard Brownian motion, $N(ds,dz)$ is an independent Poisson random measure on ${\mathbb R}_{+}\times ({\mathbb R}^{d}\setminus \{0\})$
with characteristic measure $dt \nu(dz)$. Here we assume that $\nu(dz)$ is a L\'{e}vy measure on ${\mathbb R}^{d}\setminus \{0\}$ satisfying
$\int_{{\mathbb R}^{d}\setminus \{0\}} (|z|^2\wedge 1)\nu(dz)<\infty$ with $|z|=\sqrt{\sum_{i=1}^{d} z_{i}^{2}}$.
The stochastic process $X=(X_{t}, t\geq 0)$, starting from $x_{0}
\in {\mathbb R}^{d}$, is defined as the unique strong solution to the following
stochastic differential equation (SDE)
\begin{equation}
dX_{t}=b(X_{t},\theta)dt+ \varepsilon dL_{t}, t \in [0,1]; \quad  X_{0}=x_{0}, \label{e.1.2}
\end{equation}
where $\theta \in \Theta=\bar{\Theta}_{0}$ (the closure of $\Theta_{0}$) with
$\Theta_{0}$ being an open bounded convex
subset of ${\mathbb R}^{p}$, and $b=(b_{1},\dots,b_{d}): {\mathbb R}^{d} \times \Theta \to {\mathbb R}^{d}$ is a known function. Without loss of generality, we assume that $\varepsilon \in (0,1]$. The regularity conditions
on $b$ will be provided in Section 2.
Assume that this process is observed  at regularly spaced time
points  $\{ t_{k}=k/n, \ k=1, 2, \dots, n \}$. The only unknown quantity in SDE \eref{e.1.2}
is the parameter $\theta$. Let $\theta_{0} \in \Theta_{0}$ be the true  value of the parameter $\theta$.
The purpose of this paper is to study
the least squares estimator for the true value $\theta_{0}$
based on the sampling  data $(X_{t_{k}}) _{k=1}^{n}$ with small dispersion $\varepsilon$ and large sample size $n$.

In the case of  diffusion processes driven by Brownian motion,
a popular method is the maximum likelihood estimator (MLE)
based on the Girsanov density when the processes can be observed continuously
(see Prakasa Rao \cite{PRao99},
Liptser and Shiryaev \cite{LS01}, Kutoyants \cite{Ku04}). When a diffusion process is observed only at
discrete times, in most cases the transition density and hence the likelihood function of the observations
is not explicitly computable. In order to overcome this difficulty, some approximate likelihood methods have
been proposed by Lo \cite{Lo88},
Pedersen {\cite{Ped95a}-\cite{Ped95b}, Poulsen \cite{Pou99}, and A\"{i}t-Sahalia \cite{Ait02}. For a comprehensive review on MLE and other
related methods, we refer to S{\o}rensen \cite{Sor04}.
The least squares estimator (LSE) is asymptotically equivalent to the MLE.
For the LSE,  the convergence in probability
was proved in  Dorogovcev \cite{Dor76} and
Le Breton \cite{LeB76},   the strong consistency was studied in Kasonga \cite{Kas88}, and
the asymptotic distribution was studied in Prakasa Rao \cite{PRao83}.  For a more recent comprehensive discussion, we refer to
Prakasa Rao \cite{PRao99},  Kutoyants \cite{Ku04} and  the
references therein.

The parametric estimation problems for diffusion processes with jumps based on discrete observations
 have been studied by
Shimizu and Yoshida \cite{SY06} and Shimizu \cite{Shi06} via the {quasi-}maximum likelihood. They established
consistency and asymptotic normality for the proposed estimators.
Moreover, Ogihara and Yoshida \cite{OY11} showed some stronger results than the ones by Shimizu and Yoshida \cite{SY06}, and also
investigated an adaptive Bayes-type estimator with its asymptotic properties.
The driving jump processes considered in Shimizu and Yoshida \cite{SY06}, Shimizu \cite{Shi06} and Ogihara and Yoshida \cite{OY11}
include a large class of L\'{e}vy processes such as compound Poisson processes, gamma, inverse Gaussian,
variance gamma, normal inverse Gaussian or some generalized tempered stable processes. Masuda \cite{Mas05} dealt with the
consistency
and asymptotic normality of the TFE (trajectory-fitting estimator) and LSE
 when the  driving  process is a
zero-mean adapted process (including L\'{e}vy process) with finite
moments. The parametric estimation for L\'{e}vy-driven Ornstein-Uhlenbeck processes was also studied by
Brockwell {\it et al.} \cite{BDY07}, Spiliopoulos \cite{Spi08}, and Valdivieso {\it et al.} \cite{VST09}.
However, the aforementioned papers were unable to cover an important class of driving L\'{e}vy processes, namely
$\alpha$-stable L\'{e}vy motions with $\alpha \in (0,2)$. Recently, Hu and Long \cite{HL07}-\cite{HL09} have started the
study on parameter estimation for Ornstein-Uhlenbeck processes driven by $\alpha$-stable L\'{e}vy motions. They obtained
some new asymptotic results on the proposed TFE and LSE under continuous or discrete observations, which are different from the
classical cases where asymptotic distributions are normal. Fasen \cite{Fas11} extended the results of Hu and Long \cite{HL09} to multivariate
Ornstein-Uhlenbeck processes driven by $\alpha$-stable L\'{e}vy motions. Masuda \cite{Mas10} proposed a self-weighted least absolute
deviation estimator for discretely observed ergodic Ornstein-Uhlenbeck processes driven by symmetric L\'{e}vy processes.

The asymptotic theory of parametric estimation for diffusion processes with small white noise
based on continuous-time observations has been well developed (see, e.g., Kutoyants \cite{Ku84, Ku94},
Yoshida \cite{Y92a, Y03}, Uchida and Yoshida \cite{UY04a}). There have
been many applications of small noise asymptotics to mathematical finance, see for example
Yoshida \cite{Y92b}, Takahashi \cite{Tak99}, Kunitomo and Takahashi \cite{KT01},
Takahashi and Yoshida \cite{TY04}, Uchida and Yoshida \cite{UY04b}. From a practical
point of view in parametric inference, it
is more realistic and interesting to consider asymptotic estimation for diffusion processes with
small noise based on discrete observations. Substantial progress has been  made in this direction.
Genon-Catalot \cite{GC90} and Laredo \cite{Lar90} studied the efficient estimation of drift parameters
of small diffusions from discrete observations when $\varepsilon \to 0$ and
$ n \to \infty$. S{\o}rensen \cite{Sor00} used martingale estimating
functions to establish consistency and asymptotic normality of the estimators of drift and
diffusion coefficient parameters when $\varepsilon \to 0$ and $n$ is fixed. S{\o}rensen and Uchida
\cite{SU03} and Gloter and S{\o}rensen \cite{GS09} used a contrast function to study the efficient estimation for unknown  parameters
in both drift and diffusion coefficient functions. Uchida \cite{Uch04, Uch08} used the martingale estimating function approach to
study estimation of drift parameters for small diffusions under weaker conditions. Thus, in the cases
of small diffusions, the asymptotic distributions of the estimators are normal under suitable
conditions on $\varepsilon$ and $n$.

Long \cite{Long09} studied the parameter estimation problem for discretely observed one-dimensional Ornstein-Uhlenbeck processes
with small L\'{e}vy noises. In that paper, the drift function is linear in both $x$ and $\theta$ ($(b(x,\theta)=-\theta x$),
the driving L\'{e}vy process is
$L_{t}=aB_{t}+bZ_{t}$, where $a$ and $b$ are known constants, $\{B_{t}, t\geq 0\}$ is the standard Brownian motion
and $Z_{t}$ is a $\alpha$-stable L\'{e}vy motion independent of $\{B_{t}, t\geq 0\}$.
The consistency and rate of convergence of the least squares estimator are established.
The asymptotic distribution of the LSE is shown to be the convolution of a normal distribution and a stable
distribution. In a similar framework, Long \cite{Long10} discussed the statistical estimation of the drift parameter for
a class of SDEs with special drift function $b(x,\theta)=\theta b(x)$. Ma \cite{Ma10} extended the results of Long \cite{Long09}
to the case when the driving noise is a general L\'{e}vy process. However, all the drift functions discussed in
Long \cite{Long09, Long10} and Ma \cite{Ma10} are linear in $\theta$, which restricts the applicability of their models and results.
In this paper, we allow the drift function $b(x,\theta)$ to be nonlinear in both $x$ and $\theta$, and the driving noise
to be a general L\'{e}vy process.
We are interested in estimating the drift parameter in SDE \eref{e.1.2}
based on discrete observations $\{X_{t_{i}}\}_{i=1}^{n}$ when $\varepsilon \to 0$ and $n \to \infty$.
We shall use the least squares method to obtain an asymptotically consistent estimator.

Consider the following {\it contrast function}
$$
\Psi_{n,\varepsilon}(\theta)
=\sum_{k=1}^{n}\frac{|X_{t_{k}}-X_{t_{k-1}}- b(X_{t_{k-1}}, \theta)\cdot \Delta t_{k-1}|^2}{\varepsilon^2 \Delta t_{k-1} }, $$
where $\Delta t_{k-1}=t_{k}-t_{k-1}=1/n$.
Then the LSE $\hat{\theta}_{n,\varepsilon}$ is defined as
$$
\hat{\theta}_{n,\varepsilon}:=\arg\min_{\theta\in \Theta} \Psi_{n,\varepsilon}(\theta).
$$
Since minimizing $\Psi_{n,\varepsilon} (\theta)$ is equivalent to minimizing
\[
\Phi_{n,\varepsilon}(\theta):= \varepsilon^2 (\Psi_{n,\varepsilon}(\theta)-\Psi_{n,\varepsilon}(\theta_{0})),
\]
we may write the LSE as
$$
\hat{\theta}_{n,\varepsilon}=\arg\min_{\theta \in \Theta} \Phi_{n,\varepsilon}(\theta).
$$
We shall use this fact later for convenience of the proofs.

In the nonlinear case, it is generally very difficult or impossible to obtain an explicit formula for the least squares estimator
$\hat{\theta}_{n,\varepsilon}$.  However, we can use some nice criteria in statistical inference (see Chapter 5 of Van der Vaart \cite{Van98} and
Shimizu \cite{Shi10} for a more general criterion) to establish the consistency of the LSE as well as its asymptotic behaviors (asymptotic distribution and rate of convergence). In this paper, we consider the asymptotics of the LSE $\hat{\theta}_{n,\varepsilon}$ with high frequency ($n \to \infty$)
and small dispersion ($ \varepsilon \to 0$).
Our goal is to prove that $\hat{\theta}_{n, \varepsilon} \to \theta_{0}$ in probability and to establish its rate of
convergence and asymptotic distributions. We obtain some new asymptotic distributions for
the LSE in our general setting, which are the convolutions of normal distribution and a distribution related to the jump part
of the driving L\'{e}vy process.

The paper is organized as follows.
 In Section 2, we state our main result with some remarks and examples.
We establish the consistency of the LSE $\hat{\theta}_{n, \varepsilon}$,
{ and give its asymptotic distribution, which is a natural extension of the classical small-diffusion cases}.
All the proofs are given in Section 3.

\setcounter{equation}{0}
\section{Main results}

\subsection{Notation and assumptions}

Let $X^0=(X^0_{t}, t\geq 0)$ be the solution to the underlying ordinary differential equation (ODE) under the true value of the drift parameter:
$$
dX_{t}^{0}=b(X_{t}^{0}, \theta_{0})dt, \quad  X_{0}^{0}=x_{0}.
$$
For a multi-index $m=(m_{1},\dots, m_{k})$, we define a derivative operator in $z\in \mathbb{R}^k$ as
$\partial_z^m := \partial_{z_{1}}^{m_1}\cdots \partial_{z_k}^{m_k}$,
where $\partial_{z_i}^{m_i} := \partial^{m_i}/\partial z_i^{m_i}$.
Let $C^{k,l}({\mathbb R}^{d}\times \Theta;{\mathbb R})$
be the space of all functions $f: {\mathbb R}^{d}\times \Theta \to {\mathbb R}$
which is $k$ and $l$ times continuously differentiable with respect to $x$ and $\theta$, respectively.
Moreover $C^{k,l}_\uparrow({\mathbb R}^{d}\times \Theta;{\mathbb R})$
 is a class of $f\in C^{k,l}({\mathbb R}^{d}\times \Theta;{\mathbb R})$ satisfying that
$\sup_{\theta \in \Theta}|\partial_\theta^\alpha \partial_x^\beta f(x,\theta)| \leq C(1+|x|)^{\lambda}$ for universal positive constants $C$ and $\lambda$,
where
$\alpha=(\alpha_{1},\dots, \alpha_{p})$ and $\beta=(\beta_{1},\dots,\beta_{d})$ are multi-indices with
$0\leq \sum_{i=1}^{p}\alpha_{i}\leq { l}$ and $0\leq \sum_{i=1}^{d}\beta_{i}\leq { k}$, respectively.

We introduce the following set of assumptions.
\begin{description}
\item{(A1)}
There exists a constant { $K>0$}
such that
\[
|b(x,\theta)-b(y,\theta)| \leq K|x-y|;\quad |b(x,\theta)|\leq K(1+|x|)
\]
for each $x, y \in {\mathbb R}^{d}$ and $\theta \in \Theta$.

\item{(A2)} { $b(\cdot,\cdot)\in C^{2,3}_\uparrow({\mathbb R}^{d}\times \Theta;{\mathbb R})$.}

\item{(A3)} { $\theta \ne \theta_0\ \Leftrightarrow\ b(X_{t}^{0},\theta)\ne b(X_{t}^{0},\theta_{0})$
for at least one value of $t\in [0,1]$.}

\item{(A4)} $I(\theta_{0})=(I^{ij}(\theta_{0}))_{1\leq i,j\leq p}$ is positive definite, where
$$ I^{ij}(\theta)=\int_{0}^{1}(\partial_{\theta_{i}}b)^{T}(X_{s}^{0},\theta) \partial_{\theta_{j}}b(X_{s}^{0},\theta)ds. $$
\end{description}

It is well-known that SDE (1.2) has a unique strong solution under (A1). For convenience, we shall use $C$ to denote a generic constant whose
value may vary from place to place. For a matrix $A$, we define $|A|^2={\rm tr}(AA^{T})$, where $A^{T}$ is the transpose of $A$.
In particular, $|\sigma|^2=\sum_{i=1}^{d}\sum_{j=1}^{r}\sigma_{ij}^{2}$.

\subsection{Asymptotic behavior of LSE}

{The consistency of our estimator $\hat{\theta}_{n,\varepsilon}$ is given as follows. }

\begin{Thm} \label{t.2.1}
Under conditions (A1)--(A3), we have
$$
\hat{\theta}_{n,\varepsilon} \overset{P_{\theta_{0}}}{\longrightarrow} \theta_{0}
$$
as $\varepsilon \to 0$ and $n \to \infty$.
\end{Thm}

{ The next theorem gives the asymptotic distribution of $\hat{\theta}_{n,\varepsilon}$.
As is easily seen, our result includes the case of S{\o}rensen and Uchida \cite{SU03} as a special case.}

\begin{Thm} \label{t.3.1}
Under conditions (A1)--(A4), we have
\begin{equation}\label{e.3.4}
\varepsilon^{-1}(\hat{\theta}_{n,\varepsilon}
-\theta_{0})
\overset{P_{\theta_{0}}}{\longrightarrow}
I^{-1}(\theta_{0}) S(\theta_0),
\end{equation}
as $\varepsilon \to 0$, $n \to \infty$ and $n\varepsilon \to \infty$, {where
\[
S(\theta_0):=\left(\int_{0}^{1}(\partial_{\theta_{1}}b)^{T}(X_{s}^{0},\theta_{0})dL_{s}, \dots, \int_{0}^{1}(\partial_{\theta_{p}}b)^{T}(X_{s}^{0},\theta_{0})dL_{s}\right)^{T}.
\]
}
\end{Thm}

{
\begin{remark}
One of our main contributions is that we no longer require any high-order moments condition on $X$ as in, e.g.,  S{\o}rensen and Uchida \cite{SU03} and others,
which makes our results applicable in many practical models.
\end{remark}
}

\begin{remark}\label{t.3.2}
In general, the limiting distribution on the right-hand side of \eref{e.3.4}
is a convolution of a normal distribution and a distribution related to the jump part of the L\'{e}vy process. In particular, if the driving L\'{e}vy process $L$ is the linear combination of standard Brownian motion and $\alpha$-stable motion, the limiting distribution becomes the convolution of
a normal distribution and a stable distribution.
\end{remark}

\begin{remark}
When $d=1$ and $b(x,\theta)=-\theta x$, i.e., SDE \eref{e.1.2} is linear and driven by a general L\'{e}vy process, Theorem \ref{t.3.1} reduces to Theorem 1.1 of Ma \cite{Ma10}.
When the driving L\'{e}vy process is a linear combination of standard Brownian motion and $\alpha$-stable motion, Theorem \ref{t.3.1} was discussed in Long \cite{Long09} and Ma \cite{Ma10}.
\end{remark}

\begin{remark}
Our results and arguments in the paper can be extended to the SDEs driven
by small semi-martingale noises.
\end{remark}

\begin{example}
We consider a one-dimensional stochastic process in \eref{e.1.2} with drift function $b(x,\theta)=\theta_{1}+\theta_{2}x$. We assume that
the true value $\theta_{0}=(\theta_{1}^{0},\theta_{2}^{0})$ of $\theta=(\theta_{1},\theta_{2})$ belongs to $\Theta_{0}=(c_{1}, c_{2}) \times (c_{3}, c_{4}) \subset {\mathbb R}^2$ with $c_{1}<c_{2}$ and $c_{3}<c_{4}$. Then, $X^{0}$ satisfies the following ODE \\
\begin{equation}
dX_{t}^{0}=(\theta_{1}^{0}+\theta_{2}^{0} X_{t}^{0})dt, \quad X_{0}^{0}=x_{0}. \nonumber
\end{equation}
The explicit solution is given by $X_{t}^{0}=e^{\theta_{2}^{0}t}x_{0}+\frac{\theta_{1}^{0}(e^{\theta_{2}^{0}t}-1)}{\theta_{2}^{0}}$ when $\theta_2^0\not=0$; $X_{t}^{0}=x_{0}+\theta_{1}^{0}t$ when $\theta_2^0=0$.
The LSE $\hat{\theta}_{n,\varepsilon}=(\hat{\theta}_{n,\varepsilon,1},\hat{\theta}_{n,\varepsilon,2})^{T}$ of $\theta_{0}$ is given by
\begin{eqnarray*}
&\ & \hat{\theta}_{n,\varepsilon,1}=(X_{1}-X_{0})-\hat{\theta}_{n,\varepsilon,2} \left(\frac{1}{n}\sum_{k=1}^{n}X_{t_{k-1}}\right), \\
&\ & \hat{\theta}_{n,\varepsilon,2}=\frac{\sum_{k=1}^{n}(X_{t_{k}}-X_{t_{k-1}})X_{t_{k-1}}-(X_{1}-X_{0})\left(\frac{1}{n}\sum_{k=1}^{n}X_{t_{k-1}}\right)}{\frac{1}{n}\sum_{k=1}^{n}X_{t_{k-1}}^{2} -\left(\frac{1}{n}\sum_{k=1}^{n}X_{t_{k-1}}\right)^{2}}.
\end{eqnarray*}
Note that $\partial_{\theta_{1}}b(x,\theta)=1$ and $ \partial_{\theta_{2}}b(x,\theta)=x$. In this case, the limiting random vector in Theorem 2.2 is $I^{-1}(\theta_{0})(\int_{0}^{1}dL_{s}, \int_{0}^{1}X_{s}^{0}dL_{s})^{T}$,
where
$$
I(\theta_{0})=\left(\begin{array}{cc} \int_{0}^{1}ds & \int_{0}^{1}X_{s}^{0}ds \\
\int_{0}^{1}X_{s}^{0}ds & \int_{0}^{1}(X_{s}^{0})^2 ds
\end{array} \right).
$$
\end{example}

\begin{example}
We consider a one-dimensional stochastic process in \eref{e.1.2} with drift function $b(x,\theta)=\sqrt{\theta+x^2}$. We assume that
the true value $\theta_{0}$ of $\theta$ belongs to $\Theta_{0}=(c_{1}, c_{2}) \subset {\mathbb R}$ with $0<c_{1}<c_{2}<\infty$. Then, $X^{0}$ satisfies the following ODE \\
\begin{equation}
dX_{t}^{0}=\sqrt{\theta_{0}+(X_{t}^{0})^2}dt, \quad X_{0}^{0}=x_{0}. \nonumber
\end{equation}
The explicit solution is given by $X_{t}^{0}=\frac{(x_{0}+\sqrt{\theta_{0}+x_{0}^2})^{2} e^{2t}-\theta_{0}}{2(x_{0}+\sqrt{\theta_{0}+x_{0}^2})e^{t}}$.
It is easy to verify that the LSE $\hat{\theta}_{n,\varepsilon}$ of $\theta$ is a solution to the following nonlinear equation
$$
\sum_{k=1}^{n} \frac{X_{t_{k}}-X_{t_{k-1}}}{\sqrt{\theta+X_{t_{k-1}}^2}}=1.
$$
Since it is impossible to get the explicit expression for $\hat{\theta}_{n,\varepsilon}$, we solve the above equation numerically (e.g. by
using Newton's method). Note that $\partial_{\theta} b(x,\theta)=\frac{1}{2\sqrt{\theta+x^2}}$. It is clear that the limiting random variable
in Theorem 2.2 is
$I^{-1}(\theta_{0}) \int_{0}^{1}\frac{1}{2\sqrt{\theta_{0}+(X_{s}^{0})^2}} dL_{s}$, where
$I(\theta_{0})=\int_{0}^{1}\frac{1}{4\left(\theta_{0}+(X_{s}^{0})^2\right)}ds$.
In particular, we assume that $L_{t}=aB_{t}+\sigma Z_{t}$, where $B_{t}$ is the standard Brownian motion and $Z_{t}$ is a standard $\alpha$-stable L\'{e}vy motion independent
of $B_{t}$. Let us denote by $N$ a random variable with the standard normal distribution and $U$ a random variable with the standard $\alpha$-stable
distribution $S_{\alpha}(1,\beta,0)$, where $\alpha \in (0,2)$ is the index of stability and $\beta \in [-1,1]$ is the skewness parameter.
By using the self-similarity and time change, we can easily show that the limiting random variable in Theorem 2.2 has the identical distribution
as $$a I^{-\frac{1}{2}}(\theta_{0}) N + \sigma I^{-1}(\theta_{0}) \left[\int_{0}^{1}\left(\frac{1}{2\sqrt{\theta_{0}+(X_{s}^{0})^2}}\right)^{\alpha}
ds\right]^{1/\alpha} U. $$
\end{example}

\begin{example}
We consider a two-dimensional stochastic process in \eref{e.1.2} with drift function $b(x,\theta)=C+Ax$,
where $C=(c_{1}, c_{2})^{T}$, $A=(A_{ij})_{1\leq i,j\leq 2}$ and $x=(x_{1}, x_{2})^{T}$. We assume that the eigenvalues of $A$ have positive
real parts. We want to estimate $\theta=(\theta_{1}, \dots ,\theta_{6})^{T}=(c_{1}, A_{11}, A_{12}, c_{2}, A_{21}, A_{22})^{T} \in \Theta
\subset {\mathbb R}^{6}$, whose true value is $\theta_{0}=(c_{1}^{0}, A_{11}^{0}, A_{12}^{0}, c_{2}^{0}, A_{21}^{0}, A_{22}^{0})^{T}$. Then $X_{t}^{0}$ satisfies the following ODE \\
\begin{equation}
dX_{t}^{0}=(C_{0}+A_{0}X_{t}^{0})dt, \quad X_{0}^{0}=x_{0}. \nonumber
\end{equation}
The explicit solution is given by $X_{t}^{0}=e^{A_{0}t} x_{0}+\int_{0}^{t} e^{A_{0}(t-s)}C_{0}ds$. After some basic calculation, we find that the
LSE $\hat{\theta}_{n,\varepsilon}=(\hat{\theta}_{n,\varepsilon,i})_{1\leq i\leq 6}$ is given by
$$
\left(\begin{array}{c} \hat{\theta}_{n,\varepsilon,1} \\
 \hat{\theta}_{n,\varepsilon,2} \\
  \hat{\theta}_{n,\varepsilon,3}
\end{array} \right)
= \Lambda_{n}^{-1}\left(\begin{array}{c} n\sum_{k=1}^{n} Y_{k}^{(1)} \\
 n\sum_{k=1}^{n} Y_{k}^{(1)}X_{t_{k-1}}^{(1)} \\
  n\sum_{k=1}^{n} Y_{k}^{(1)}X_{t_{k-1}}^{(2)}
\end{array} \right)
  \mbox{\rm and} \
\left(\begin{array}{c} \hat{\theta}_{n,\varepsilon,4} \\
 \hat{\theta}_{n,\varepsilon,5} \\
  \hat{\theta}_{n,\varepsilon,6}
\end{array} \right)
= \Lambda_{n}^{-1}\left(\begin{array}{c} n\sum_{k=1}^{n} Y_{k}^{(2)} \\
 n\sum_{k=1}^{n} Y_{k}^{(2)}X_{t_{k-1}}^{(1)} \\
  n\sum_{k=1}^{n} Y_{k}^{(2)}X_{t_{k-1}}^{(2)}
\end{array} \right),
$$
where $X_{t_{k-1}}^{(i)}$ ($i=1,2$) are the components of $X_{t_{k-1}}$, $Y_{k}^{(i)}$ ($i=1,2$) are the components of
$Y_{k}=X_{t_{k}}-X_{t_{k-1}}$, and
$$
\Lambda_{n}=\left(\begin{array}{ccc} n & \sum_{k=1}^{n}X_{t_{k-1}}^{(1)} & \sum_{k=1}^{n}X_{t_{k-1}}^{(2)} \\
\sum_{k=1}^{n}X_{t_{k-1}}^{(1)} & \sum_{k=1}^{n}\left(X_{t_{k-1}}^{(1)}\right)^2 & \sum_{k=1}^{n}X_{t_{k-1}}^{(1)}X_{t_{k-1}}^{(2)} \\
\sum_{k=1}^{n}X_{t_{k-1}}^{(2)} & \sum_{k=1}^{n}X_{t_{k-1}}^{(1)}X_{t_{k-1}}^{(2)} &\sum_{k=1}^{n}\left(X_{t_{k-1}}^{(2)}\right)^2
\end{array} \right).
$$
Since it is easy and straightforward to compute the partial derivatives $\partial_{\theta_{i}}b(x,\theta)$, $1\leq i\leq 6$, and the limiting
random vector in Theorem 2.2, we omit the details here.
\end{example}

\setcounter{equation}{0}
\section{Proofs}

\subsection{Proof of Theorem \ref{t.2.1}}

We first establish some preliminary lemmas.
In the sequel, we shall use the notation
\[
Y_{t}^{n,\varepsilon}:=X_{[nt]/n}
\]
for the stochastic process $X$ defined by \eref{e.1.2}, where $[nt]$ denotes the integer part of $nt$.

\begin{lemma}\label{t.2.2}
The sequence $\{Y_{t}^{n,\varepsilon}\}$ converges to the deterministic process $\{X_{t}^{0}\}$ uniformly on compacts in probability as $\varepsilon \to 0$ and $n \to \infty$.
\end{lemma}
{\bf Proof}. Note that
\begin{equation}\label{e.2.3}
X_{t}-X_{t}^{0}=\int_{0}^{t}(b(X_{s},\theta_{0})-b(X_{s}^{0},\theta_{0}))ds+\varepsilon L_{t}.
\end{equation}
By the Lipschitz condition on $b(\cdot)$ in (A1) and the Cauchy-Schwarz inequality, we find that
\begin{eqnarray*}
|X_{t}-X_{t}^{0}|^2 &\leq & 2\left|\int_{0}^{t}(b(X_{s},\theta_{0})-b(X_{s}^{0},\theta_{0}))ds\right|^{2}+2\varepsilon^2 |L_{t}|^2 \\
&\leq & 2t \int_{0}^{t}|b(X_{s},\theta_{0})-b(X_{s}^{0},\theta_{0})|^2ds+2\varepsilon^2 \sup_{0\leq s \leq t}|L_{s}|^2 \\
&\leq & 2K^2t\int_{0}^{t}|X_{s}-X_{s}^{0}|^2ds+2\varepsilon^2 \sup_{0\leq s \leq t}|L_{s}|^2 .
\end{eqnarray*}
By Gronwall's inequality, it follows that
$$
|X_{t}-X_{t}^{0}|^2 \leq 2\varepsilon^{2} e^{2K^{2}t^{2}} \sup_{0\leq s \leq t}|L_{s}|^{2}
$$
and consequently
\begin{equation}\label{e.2.5}
\sup_{ 0\leq t \leq T}|X_{t}-X_{t}^{0}| \leq \sqrt{2}\varepsilon e^{K^{2}T^{2}} \sup_{0\leq t \leq T}|L_{t}|,
\end{equation}
which goes to zero in probability as $\varepsilon \to 0$ for each $T>0$. Since $[nt]/n \to t$ as $n \to \infty$, we conclude that the statement holds.\\

\smallskip\noindent
\begin{lemma}\label{t.2.3}
Let $\tau_{m}^{n,\varepsilon}=\inf \{t\geq 0: |X_{t}^{0}|\geq m \ \mbox{\rm or} \  |Y_{t}^{n,\varepsilon}|\geq m\}$. Then,
$\tau_{m}^{n,\varepsilon} \to \infty\ {a.s.}$ uniformly in $n$ and $\varepsilon$ as $m \to \infty$.
\end{lemma}
{\bf Proof}. Note that
$$X_{t}=x_{0}+\int_{0}^{t} b(X_{s},\theta_{0})ds+\varepsilon L_{t}. $$
By the linear growth condition on $b$ and the Cauchy-Schwarz inequality, we get
\begin{eqnarray*}
|X_{t}|^{2} &\leq & 2(|x_{0}|+\varepsilon |L_{t}|)^{2}+2\left|\int_{0}^{t} b(X_{s},\theta_{0})ds\right|^{2} \\
&\leq & 2\left(|x_{0}|+\varepsilon \sup_{0\leq s \leq t}|L_{s}|\right)^{2}+2t\int_{0}^{t} |b(X_{s},\theta_{0})|^2ds \\
&\leq & 2\left(|x_{0}|+\varepsilon \sup_{0\leq s \leq t}|L_{s}|\right)^{2}+2K^{2}t\int_{0}^{t} (1+|X_{s}|)^2ds \\
&\leq & \left[2(|x_{0}|+\varepsilon \sup_{0\leq s \leq t}|L_{s}|)^{2}+4K^2t^{2}\right]+4K^{2}t\int_{0}^{t}|X_{s}|^2ds.
\end{eqnarray*}
Gronwall's inequality yields that
$$|X_{t}|^2 \leq \left[2(|x_{0}|+\varepsilon \sup_{0\leq s \leq t}|L_{s}|)^{2}+4K^2t^{2}\right]e^{4K^{2}t^2} $$
and
$$ |X_{t}|
 \leq \left[\sqrt{2}(|x_{0}|+\varepsilon \sup_{0\leq s \leq t}|L_{s}|)+2Kt\right]e^{2K^{2}t^2}. $$
Thus, it follows that
$$ |Y_{t}^{n,\varepsilon}|=|X_{[nt]/n}| \leq
\left[\sqrt{2}(|x_{0}|+\sup_{0\leq s \leq t}|L_{s}|)+2Kt\right]e^{2K^{2}t^2},
$$
which is almost surely finite. Therefore the proof is complete. $\hfill \square $ \\

\smallskip\noindent
We shall use $\nabla_{x}f(x,\theta)=(\partial_{x_{1}}f(x,\theta), \dots, \partial_{x_{d}}f(x,\theta))^{T}$ to denote the gradient
operator of $f(x,\theta)$ with respect to $x$.\\

\smallskip\noindent
\begin{lemma}\label{t.2.4}
Let $f \in C^{1,1}_{\uparrow}({\mathbb R}^{d}\times \Theta;{\mathbb R})$. Assume (A1)-(A2). Then, we have
$$ \frac{1}{n} \sum_{k=1}^{n} f(X_{t_{k-1}},\theta) \overset{P_{\theta_{0}}}{\longrightarrow} \int_{0}^{1} f(X_{s}^{0},\theta)ds $$
as $\varepsilon \to 0$ and $n \to \infty$, uniformly in $\theta \in \Theta$.
\end{lemma}
{\bf Proof}.
By the differentiability of the function $f(x,\theta)$ and Lemma 3.1, we find that
\begin{align*}
\sup_{\theta \in \Theta}&\left|\frac{1}{n}\sum_{k=1}^{n}f(X_{t_{k-1}},\theta)-\int_{0}^{1}f(X_{s}^{0},\theta)ds\right| \\
&=\sup_{\theta \in \Theta}\left|\int_{0}^{1}f(Y_{s}^{n,\varepsilon},\theta)ds-\int_{0}^{1}f(X_{s}^{0},\theta)ds\right| \\
&\leq \sup_{\theta \in \Theta} \int_{0}^{1} |f(Y_{s}^{n,\varepsilon},\theta)-f(X_{s}^{0},\theta)|ds \\
&\leq \sup_{\theta \in \Theta} \int_{0}^{1} \left|\int_{0}^{1}(\nabla_{x}f)^{T}(X_{s}^{0}+u(Y_{s}^{n,\varepsilon}-X_{s}^{0}),\theta)
\cdot (Y_{s}^{n,\varepsilon}-X_{s}^{0})du\right| ds \\
&\leq \int_{0}^{1} \left(\int_{0}^{1}\sup_{\theta \in \Theta}|\nabla_{x}f(X_{s}^{0}+u(Y_{s}^{n,\varepsilon}-X_{s}^{0}),\theta)|
du \right) |Y_{s}^{n,\varepsilon}-X_{s}^{0}|ds \\
&\leq \int_{0}^{1} C(1+|X_{s}^{0}|+|Y_{s}^{n,\varepsilon}|)^{\lambda}|Y_{s}^{n,\varepsilon}-X_{s}^{0}|ds \\
&\leq C\left(1+\sup_{0\le s\leq 1}|X_{s}^{0}|+\sup_{0\leq s \leq 1}|X_{s}|\right)^{\lambda} \sup_{0\leq s \leq 1}|Y_{s}^{n,\varepsilon}-X_{s}^{0}| \\
&\overset{P_{\theta_{0}}}{\longrightarrow}  0
\end{align*}
as $\varepsilon \to 0$ and $n \to \infty$. $\hfill \square $ \\

\smallskip\noindent
\begin{lemma}\label{t.2.5}
Let $f \in C^{1,1}_{\uparrow}({\mathbb R}^{d}\times \Theta;{\mathbb R})$. Assume (A1)-(A2). Then, we have that for each $ 1\leq i \leq d$ and each $\theta \in \Theta$,
$$ \sum_{k=1}^{n} f(X_{t_{k-1}},\theta)(L_{t_{k}}^{i}-L_{t_{k-1}}^{i}) \overset{P_{\theta_{0}}}{\longrightarrow} \int_{0}^{1} f(X_{s}^{0},\theta)dL_{s}^{i} $$
as $\varepsilon \to 0$ and $n \to \infty$, where
$$ L_{t}^{i}=a_{i}t+\sum_{j=1}^{r}\sigma_{ij}B_{t}^{j}+\int_{0}^{t}\int_{|z|\leq 1}z_{i}\tilde{N}(ds,dz)
+\int_{0}^{t}\int_{|z|>1}z_{i}N(ds,dz) $$
is the $i$-th component of $L_{t}$.
\end{lemma}
{\bf Proof}. Note that $$\sum_{k=1}^{n} f(X_{t_{k-1}},\theta)(L^{i}_{t_{k}}-L^{i}_{t_{k-1}})=\int_{0}^{1}f(Y_{s}^{n,\varepsilon},\theta)dL^{i}_{s}.
$$
Let $\tilde{L}^{i}_{t}=L^{i}_{t}-\int_{0}^{t}\int_{|z|>1} z_{i}N(ds,dz)$. Then, we have the following decomposition
\begin{align*}
\int_{0}^{1}f(Y_{s}^{n,\varepsilon},\theta)dL^{i}_{s} -\int_{0}^{1} f(X_{s}^{0},\theta)dL^{i}_{s} 
&= \int_{0}^{1}\int_{|z|>1} (f(Y_{s}^{n,\varepsilon},\theta)-f(X_{s}^{0},\theta)) z_{i}N(ds,dz) \nonumber \\
&\quad  +\int_{0}^{1}(f(Y_{s}^{n,\varepsilon},\theta)-f(X_{s}^{0},\theta))d\tilde{L}^{i}_{s}.
\end{align*}
Similar to the proof of Lemma 3.3, we have
\begin{align*}
\Bigg|\int_{0}^{1}&\int_{|z|>1} (f(Y_{s}^{n,\varepsilon},\theta)-f(X_{s}^{0},\theta)) z_{i}N(ds,dz) \Bigg| \\
&\leq  \int_{0}^{1}\int_{|z|>1} |f(Y_{s}^{n,\varepsilon},\theta)-f(X_{s}^{0},\theta)| |z_{i}|N(ds,dz) \\
&\leq  \int_{0}^{1}\int_{|z|>1} C(1+|X_{s}^{0}|+|Y_{s}^{n,\varepsilon}|)^{\lambda} |Y_{s}^{n,\varepsilon}-X_{s}^{0}| |z_{i}|N(ds,dz) \\
&\leq  C\left(1+\sup_{0\leq s\leq 1}|X_{s}^{0}|+\sup_{0\leq s \leq 1} |X_{s}|\right)^{\lambda} \sup_{0\leq s \leq 1}|Y_{s}^{n,\varepsilon}-X_{s}^{0}|
\int_{0}^{1}\int_{|z|>1}|z_{i}|N(ds,dz),
\end{align*}
which converges to zero in probability as $\varepsilon \to 0$ and $n \to \infty$ by Lemma 3.1. By using the stopping time $\tau_{m}^{n,\varepsilon}$, Lemma 3.1, Markov inequality and dominated convergence, we find that for any given $\eta>0$ and some fixed $m$
\begin{align}
P&\left(\left|\int_{0}^{1}(f(Y_{s}^{n,\varepsilon},\theta)-f(X_{s}^{0},\theta))1_{\{s\leq \tau_{m}^{n,\varepsilon}\}}d\tilde{L}^{i}_{s}\right|
> \eta\right) \nonumber \\
&\leq  \frac{|a_{i}|}{\eta} \int_{0}^{1}{\mathbb E}\left[|f(Y_{s}^{n,\varepsilon},\theta)-f(X_{s}^{0},\theta)|1_{\{s\leq \tau_{m}^{n,\varepsilon}\}}\right]ds \nonumber \\
&\quad + \frac{\sqrt{\sum_{j=1}^{r}\sigma_{ij}^2}}{\eta}\left(\int_{0}^{1}{\mathbb E}\left[|f(Y_{s}^{n,\varepsilon},\theta)-f(X_{s}^{0},\theta)|^{2}1_{\{s\leq \tau_{m}^{n,\varepsilon}\}}\right]ds \right)^{1/2} \nonumber \\
&\quad +\frac{1}{\eta}\left(\int_{0}^{1}{\mathbb E}\left[|f(Y_{s}^{n,\varepsilon},\theta)-f(X_{s}^{0},\theta)|^{2}1_{\{s\leq \tau_{m}^{n,\varepsilon}\}}\right]ds \cdot \int_{|z|\leq 1} |z_{i}|^2 \nu(dz) \right)^{1/2}, \label{e.2.7}
\end{align}
which goes to zero as $\vare \to 0$ and $n \to \infty$.
Then, we have
\begin{align*}
P&\left(\left|\int_{0}^{1}(f(Y_{s}^{n,\varepsilon},\theta)-f(X_{s}^{0},\theta))d\tilde{L}_{s}^{i}\right|
> \eta\right) \\
&\leq  P(\tau_{m}^{n,\varepsilon}<1)+
P\left(\left|\int_{0}^{1}(f(Y_{s}^{n,\varepsilon},\theta)-f(X_{s}^{0},\theta))1_{\{s\leq \tau_{m}^{n,\varepsilon}\}}d\tilde{L}_{s}^{i}\right|
> \eta\right),
\end{align*}
which converges to zero as $\varepsilon \to 0$ and $n \to \infty$ by Lemma 3.2 and \eref{e.2.7}. This completes the proof.
$\hfill \square$ \\

\smallskip\noindent
\begin{lemma}\label{t.2.7}
Let $f\in C_{\uparrow}^{1,1}(({\mathbb R}^{d}\times \Theta;{\mathbb R})$. Assume (A1)-(A2). Then, we have that for $1\leq i \leq d$,
$$ \sum_{k=1}^{n} f(X_{t_{k-1}},\theta)(X^{i}_{t_{k}}-X^{i}_{t_{k-1}}-b_{i}(X_{t_{k-1}},\theta_{0})\Delta t_{k-1}) \overset{P_{\theta_{0}}}{\longrightarrow} 0 $$
as $\varepsilon \to 0$ and $n \to \infty$, uniformly in $\theta \in \Theta$, where $X_{t}^{i}$ and $b_{i}$ are the $i$-th components of
$X_{t}$ and $b$, respectively.\\
\end{lemma}
{\bf Proof}. Note that
$$ X_{t_k}^{i}=X_{t_{k-1}}^{i}+\int_{t_{k-1}}^{t_{k}}b_{i}(X_{s},\theta_{0})ds+\varepsilon (L_{t_{k}}^{i}-L_{t_{k-1}}^{i}). $$
It is easy to see that
\begin{align*}
\sum_{k=1}^{n} & f(X_{t_{k-1}},\theta)(X_{t_{k}}^{i}-X_{t_{k-1}}^{i}-b_{i}(X_{t_{k-1}},\theta_{0})\Delta t_{k-1}) \nonumber \\
&= \sum_{k=1}^{n}\int_{t_{k-1}}^{t_{k}}f(X_{t_{k-1}},\theta)(b_{i}(X_{s},\theta_{0})-b_{i}(X_{t_{k-1}},\theta_{0}))ds \nonumber \\
&\quad + \vare \sum_{k=1}^{n} f(X_{t_{k-1}},\theta)(L_{t_{k}}^{i}-L_{t_{k-1}}^{i}) \nonumber \\
&= \int_{0}^{1}f(Y_{s}^{n,\varepsilon},\theta)(b_{i}(X_{s},\theta_{0})-b_{i}(Y_{s}^{n,\varepsilon},\theta_{0}))ds
+\varepsilon \int_{0}^{1}f(Y_{s}^{n,\varepsilon},\theta)dL_{s}^{i}.
\end{align*}
By the given condition on $f$ and the Lipschitz condition on $b$, we have
\begin{align*}
\sup_{\theta \in \Theta}&\left|\int_{0}^{1}f(Y_{s}^{n,\varepsilon},\theta)(b_{i}(X_{s},\theta_{0})-b_{i}(Y_{s}^{n,\varepsilon},\theta_{0}))ds\right|
\nonumber \\
&\leq  \int_{0}^{1} \sup_{\theta \in \Theta}|f(Y_{s}^{n,\varepsilon},\theta)| \cdot K|X_{s}-Y_{s}^{n,\varepsilon}|ds
\nonumber \\
&\leq  K C\int_{0}^{1}(1+|Y_{s}^{n,\varepsilon}|)^{\lambda}(|X_{s}-X_{s}^{0}|+|Y_{s}^{n,\varepsilon}-X_{s}^{0}|)ds \nonumber \\
&\leq  K C\left(1+\sup_{0\leq t\leq 1}|X_{t}|\right)^{\lambda}(\sup_{0\leq s\leq 1}|X_{s}-X_{s}^{0}|+\sup_{0\leq s\leq 1}|Y_{s}^{n,\varepsilon}-X_{s}^{0}|),
\end{align*}
which converges to zero in probability as  $\varepsilon \to 0$ and $n \to \infty$ by Lemma 3.1.
Next using the decomposition of $L_{t}$, we have
\begin{align*}
\sup_{\theta \in \Theta}&\left|\varepsilon \int_{0}^{1}f(Y_{s}^{n,\varepsilon},\theta)dL_{s}^{i}\right| \nonumber \\
&\leq \varepsilon \sup_{\theta \in \Theta}\left|a_{i}\int_{0}^{1}f(Y_{s}^{n,\varepsilon},\theta)ds\right|
+\varepsilon\sup_{\theta \in \Theta}\left|\int_{0}^{1}f(Y_{s}^{n,\varepsilon},\theta)\sum_{j=1}^{r}\sigma_{ij} dB_{s}^{j}\right| \nonumber \\
&\quad +\varepsilon\sup_{\theta \in \Theta}\left|\int_{0}^{1}\int_{|z|\leq 1}f(Y_{s}^{n,\varepsilon},\theta)z_{i}\tilde{N}(ds,dz)\right|
\nonumber \\
&\quad + \varepsilon \sup_{\theta \in \Theta}\left|\int_{0}^{1}\int_{|z|>1}f(Y_{s}^{n,\varepsilon},\theta)z_{i}N(ds,dz)\right|.
\end{align*}
It is clear that
\begin{eqnarray*}
\varepsilon \sup_{\theta \in \Theta}\left|a_{i}\int_{0}^{1}f(Y_{s}^{n,\varepsilon},\theta)ds\right|
&\leq & \varepsilon |a_{i}| C\int_{0}^{1} (1+|Y_{s}^{n,\varepsilon}|)^{\lambda}ds \nonumber \\
&\leq & \varepsilon |a_{i}| C\left(1+\sup_{0\leq s \leq 1}|X_{s}|\right)^{\lambda},
\end{eqnarray*}
which converges to zero in probability as $\varepsilon \to 0$ and $n \to \infty$, and
\begin{align*}
\varepsilon \sup_{\theta \in \Theta}&\left|\int_{0}^{1}\int_{|z|>1}f(Y_{s}^{n,\varepsilon},\theta)z_{i}N(ds,dz)\right| \nonumber \\
&\leq  \varepsilon \int_{0}^{1}\int_{|z|>1}\sup_{\theta \in \Theta}|f(Y_{s}^{n,\varepsilon},\theta)| \cdot |z_{i}|N(ds,dz) \nonumber \\
&\leq  \varepsilon \int_{0}^{1}\int_{|z|>1}C (1+|Y_{s}^{n,\varepsilon}|)^{\lambda}\cdot |z_{i}|N(ds,dz) \nonumber \\
&\leq  \varepsilon C \left(1+\sup_{0\leq s \leq 1}|X_{s}|\right)^{\lambda} \int_{0}^{1}\int_{|z|>1}|z_{i}|N(ds,dz),
\end{align*}
which converges to zero in probability. Note that
\begin{align}
P&\left(\varepsilon\sup_{\theta \in \Theta}\left|\int_{0}^{1}f(Y_{s}^{n,\varepsilon},\theta)\sum_{j=1}^{r}\sigma_{ij} dB_{s}^{j}\right|>\eta \right) \nonumber \\
&\leq P(\tau_{m}^{n,\varepsilon}<1) +
P\left(\varepsilon\sup_{\theta \in \Theta}\left|\int_{0}^{1}f(Y_{s}^{n,\varepsilon},\theta)1_{\{s\leq \tau_{m}^{n,\varepsilon}\}}\sum_{j=1}^{r}\sigma_{ij} dB_{s}^{j}\right|>\eta\right).
\label{e.2.14}
\end{align}
Let $$u_{n,\varepsilon}^{i}(\theta)=\varepsilon
\int_{0}^{1}f(Y_{s}^{n,\varepsilon},\theta)1_{\{s\leq
\tau_{m}^{n,\varepsilon}\}}\sum_{j=1}^{r}\sigma_{ij} dB_{s}^{j}, \
1\leq i \leq d.$$ We want to prove that
$u_{n,\varepsilon}^{i}(\theta) \to 0$ in probability as
$\varepsilon \to 0$ and $n \to \infty$, uniformly in $\theta \in
\Theta$. It suffices to show the pointwise convergence and the
tightness of the sequence $\{u_{n,\varepsilon}^{i}(\cdot)\}$. For
the pointwise convergence, by the Chebyshev inequality and Ito's
isometry, we have
\begin{eqnarray}
&\ & P(|u_{n,\varepsilon}^{i}(\theta)|>\eta) \nonumber \\
&\leq & \varepsilon^2 \eta^{-2} {\mathbb E}\left[\left|\int_{0}^{1}f(Y_{s}^{n,\varepsilon},\theta)1_{\{s\leq \tau_{m}^{n,\varepsilon}\}}\sum_{j=1}^{r}\sigma_{ij} dB_{s}^{j}\right|^2\right] \nonumber \\
&\leq & \left(\sum_{j=1}^{r}\sigma_{ij}^{2}\right)\varepsilon^2
\eta^{-2} \int_{0}^{1}{\mathbb
E}\left[|f(Y_{s}^{n,\varepsilon},\theta)|^2
1_{\{s\leq \tau_{m}^{n,\varepsilon}\}}\right]ds \nonumber \\
&\leq & \left(\sum_{j=1}^{r}\sigma_{ij}^{2}\right)\varepsilon^2
\eta^{-2} \int_{0}^{1}{\mathbb
E}\left[C^2(1+|Y_{s}^{n,\varepsilon}|)^{2\lambda}
1_{\{s\leq \tau_{m}^{n,\varepsilon}\}}\right]ds \nonumber \\
&\leq & \left(\sum_{j=1}^{r}\sigma_{ij}^{2}\right)\varepsilon^2
\eta^{-2}C^2 (1+m)^{2\lambda}, \label{e.2.15}
\end{eqnarray}
which converges to zero as $\varepsilon \to 0$ and $n \to \infty$ with fixed $m$. For the tightness of $\{u_{n,\varepsilon}^{i}(\cdot)\}$, by using Theorem 20 in Appendix I of Ibragimov and Has'minskii \cite{IH81}, it is enough to prove the following two inequalities
\begin{eqnarray}
&\ & {\mathbb E}[|u_{n,\varepsilon}^{i} (\theta)|^{2q}] \leq C, \label{e.2.16} \\
&\ & {\mathbb E}[|u_{n,\varepsilon}^{i} (\theta_{2})-u_{n,\varepsilon}^{i} (\theta_{1})|^{2q}]\leq C|\theta_{2}-\theta_{1}|^{2q}
\label{e.2.17}
\end{eqnarray}
for $\theta, \theta_{1},\theta_{2} \in \Theta$, where $2q>p$. The
proof of \eref{e.2.16} is very similar to moment estimates in
\eref{e.2.15} by replacing Ito's isometry with the
Burkholder-Davis -Gundy inequality. So we omit the details here.
For \eref{e.2.17}, by using Taylor's formula and the
Burkholder-Davis-Gundy inequality, we have
\begin{eqnarray*}
&\ & {\mathbb E}[|u_{n,\varepsilon}^{i} (\theta_{2})-u_{n,\varepsilon}^{i} (\theta_{1})|^{2q}] \\
&\leq & \varepsilon^{2q} C_{q}
\left(\sum_{j=1}^{r}\sigma_{ij}^2\right)^{q}{\mathbb
E}\left[\left(\int_{0}^{1}(f(Y_{s}^{n,\varepsilon},\theta_{2})
-f(Y_{s}^{n,\varepsilon},\theta_{1}))^{2}1_{\{s\leq \tau_{m}^{n,\varepsilon}\}}ds\right)^{q}\right] \\
&\leq & \varepsilon^{2q} C_{q}
\left(\sum_{j=1}^{r}\sigma_{ij}^2\right)^{q}{\mathbb
E}\left[\left(\int_{0}^{1}\int_{0}^{1}|\theta_{2}-\theta_{1}|^{2}
|\nabla_{\theta}f(Y_{s}^{n,\varepsilon},\theta_{1}+v(\theta_{2}-\theta_{1}))|^{2}1_{\{s\leq
\tau_{m}^{n,\varepsilon}\}}dvds\right)^{q}\right]
\\
&\leq & \varepsilon^{2q} C_{q}
\left(\sum_{j=1}^{r}\sigma_{ij}^2\right)^{q}C^{2q}|\theta_{2}-\theta_{1}|^{2q}
{\mathbb
E}\left[\left(\int_{0}^{1}(1+|Y_{s}^{n,\varepsilon}|)^{2\lambda}
1_{\{s\leq \tau_{m}^{n,\varepsilon}\}}ds\right)^{q}\right] \\
&\leq & \varepsilon^{2q} C_{q}
\left(\sum_{j=1}^{r}\sigma_{ij}^2\right)^{q}C^{2q} (1+m)^{2\lambda
q}|\theta_{2}-\theta_{1}|^{2q}.
\end{eqnarray*}
Combining \eref{e.2.14} and the above arguments, we have that $\varepsilon\sup_{\theta \in \Theta}\left|\int_{0}^{1}f(Y_{s}^{n,\varepsilon},\theta)\sum_{j=1}^{r}\sigma_{ij} dB_{s}^{j}\right|$ converges to zero in probability
as $\varepsilon \to 0$ and $n \to \infty$. Similarly, we can prove that
$\varepsilon\sup_{\theta \in \Theta}\left|\int_{0}^{1}\int_{|z|\leq 1}f(Y_{s}^{n,\varepsilon},\theta)z_{i}\tilde{N}(ds,dz)\right|$ converges  to zero in probability
as $\varepsilon \to 0$ and $n \to \infty$.
Therefore, the proof is complete. $\hfill \square $ \\

\smallskip\noindent
Now we are  in a position to prove Theorem 2.1.\\

\smallskip\noindent
{\bf Proof of Theorem 2.1}. Note that
\begin{eqnarray*}
\Phi_{n,\varepsilon}(\theta))&=& -2\sum_{k=1}^{n}(b(X_{t_{k-1}},\theta)-
b(X_{t_{k-1}},\theta_{0}))^{T}(X_{t_{k}}-X_{t_{k-1}}-n^{-1}b(X_{t_{k-1}},\theta_{0})) \nonumber \\
&\ & +\frac{1}{n}\sum_{k=1}^{n}|b(X_{t_{k-1}},\theta)-b(X_{t_{k-1}},\theta_{0})|^2 . \nonumber \\
&:=& \Phi_{n,\varepsilon}^{(1)}(\theta)+\Phi_{n,\varepsilon}^{(2)}(\theta).
\end{eqnarray*}
By Lemma 3.5 and let $f(x,\theta)=b_{i}(x,\theta)-b_{i}(x,\theta_{0})$ ($1\leq i \leq d$), we have
$ \sup_{\theta \in \Theta}|\Phi_{n,\varepsilon}^{(1)}(\theta)| \overset{P_{\theta_{0}}}{\longrightarrow} 0$ as $\varepsilon \to 0$ and $n \to \infty$. By using Lemma 3.3
with $f(x,\theta)=|b(x,\theta)-b(x,\theta_{0})|^2$, we find
$\sup_{\theta \in \Theta}|\Phi_{n,\varepsilon}^{(2)}(\theta)-F(\theta)| \overset{P_{\theta_{0}}}{\longrightarrow} 0$ as $\varepsilon \to 0$ and $n \to \infty$,
where $F(\theta)=\int_{0}^{1}|b(X_{t}^{0},\theta)-b(X_{t}^{0},\theta_{0})|^2 dt$.
Thus combining the previous arguments, we have
$$\sup_{\theta \in \Theta}|\Phi_{n,\varepsilon}(\theta)-F(\theta)| \overset{P_{\theta_{0}}}{\longrightarrow} 0$$ as $\varepsilon \to 0$ and $n \to \infty$,
and that (A3) and the continuity of $X^0$ yield that
\[
\inf_{|\th-\th_0|>\delta}F(\theta) > F(\th_0)=0,
\]
for each $\delta>0$.
Therefore, by Theorem 5.9 of van der Vaart \cite{Van98}, we have the desired consistency, i.e., $\hat{\theta}_{n,\varepsilon} \overset{P_{\theta_{0}}}{\longrightarrow} {\theta_0}$ as
$\varepsilon \to 0$ and $n \to \infty$. This completes the proof. $\hfill \square $

\subsection{Proof of Theorem \ref{t.3.1}}

Note that
$$\nabla_{\theta}\Phi_{n,\varepsilon}(\theta)=
-2\sum_{k=1}^{n}(\nabla_{\theta}
b)^{T}(X_{t_{k-1}},\theta)(X_{t_{k}}-X_{t_{k-1}}-b(X_{t_{k-1}},\theta)\Delta t_{k-1}).$$
Let $G_{n,\varepsilon}(\theta)=(G_{n,\varepsilon}^{1},\dots, G_{n,\varepsilon}^{p})^{T}$ with
$$
G_{n,\varepsilon}^{i}(\theta)=\sum_{k=1}^{n}(\partial_{\theta_{i}}
b)^{T}(X_{t_{k-1}},\theta)(X_{t_{k}}-X_{t_{k-1}}-b(X_{t_{k-1}},\theta)\Delta t_{k-1}), \ i=1,\dots, p,
$$
and let $K_{n,\varepsilon}(\theta)=\nabla_{\theta}G_{n,\varepsilon}(\theta)$, which is a $p\times p$ matrix consisting of elements
$K_{n,\varepsilon}^{ij}(\theta)=\partial_{\theta_{j}}G_{n,\varepsilon}^{i}(\theta), 1\leq i, j \leq p$.
Moreover, we introduce the following function
$$ K^{ij}(\theta) =\int_{0}^{1}(\partial_{\theta_{j}}\partial_{\theta_{i}}b)^{T}(X_{s}^{0},\theta)(b(X_{s}^{0},\theta_{0})
-b(X_{s}^{0},\theta))ds-I^{ij}(\theta), \ 1\leq i, j \leq p. $$
Then we define the matrix function $K(\theta)=(K^{ij}(\theta))_{1\leq i,j\leq p}$.

Before proving Theorem \ref{t.3.1}, we prepare some preliminary results.

\begin{lemma}\label{t.3.3}
Assume (A1)-(A2). Then, we have that for each $i=1,\dots,p$
$$ \varepsilon^{-1} G_{n,\varepsilon}^{i}(\theta_{0}) \overset{P_{\theta_{0}}}{\longrightarrow}
\int_{0}^{1}(\partial_{\theta_{i}}b)^{T}(X_{s}^{0},\theta_{0})dL_{s} $$
as $\varepsilon \to 0$, $n \to \infty$ and $n\varepsilon \to \infty$. \\
\end{lemma}
{\bf Proof}. Note that for $1\leq i \leq p$
\begin{eqnarray*}
\varepsilon^{-1} G_{n,\varepsilon}^{i}(\theta_{0})&=&
\varepsilon^{-1}\sum_{k=1}^{n}(\partial_{\theta_{i}}
b)^{T}(X_{t_{k-1}},\theta_{0})(X_{t_{k}}-X_{t_{k-1}}-b(X_{t_{k-1}},\theta_{0})\Delta t_{k-1}) \nonumber \\
&=& \varepsilon^{-1}\sum_{k=1}^{n}(\partial_{\theta_{i}}b)^{T}(X_{t_{k-1}},\theta_{0})\int_{t_{k-1}}^{t_{k}}(b(X_{s},\theta_{0})-b(X_{t_{k-1}},\theta_{0}))ds \nonumber \\
&\ & +\sum_{k=1}^{n}(\partial_{\theta_{i}}b)^{T}(X_{t_{k-1}},\theta_{0})(L_{t_{k}}-L_{t_{k-1}}) \nonumber \\
&:=& H_{n,\varepsilon}^{(1)}(\theta_{0})+H_{n,\varepsilon}^{(2)}(\theta_{0}).
\end{eqnarray*}
By using Lemma 3.4 and letting $f(x,\theta)=\partial_{\theta_{i}}b_{j}(x,\theta)$ ($1\leq i \leq p, \ 1\leq j \leq d$) with $\theta=\theta_{0}$, we have
$$
H_{n,\varepsilon}^{(2)}(\theta_{0})=\int_{0}^{1}(\partial_{\theta_{i}}b)^{T}(Y_{s}^{n,\varepsilon},\theta_{0})dL_{s}
\overset{P_{\theta_{0}}}{\longrightarrow} \int_{0}^{1}(\partial_{\theta_{i}}b)^{T}(X_{s}^{0},\theta_{0})dL_{s}
$$
as $\varepsilon \to 0$ and $n \to \infty$. It suffices to prove that $H_{n,\varepsilon}^{(1)}(\theta_{0})$
converges to zero in probability. For $H_{n,\varepsilon}^{(1)}(\theta_{0})$, we need some delicate estimate
for the process $X_{t}$. For $s \in [t_{k-1},t_{k}]$, we have
$$ X_{s}-X_{t_{k-1}}=\int_{t_{k-1}}^{s}(b(X_{u},\theta_{0})-b(X_{t_{k-1}},\theta_{0}))du
+b(X_{t_{k-1}},\theta_{0})(s-t_{k-1})+\varepsilon (L_{s}-L_{t_{k-1}}). $$
By the Lipschitz condition on $b$ and the Cauchy-Schwarz inequality, we find that
\begin{eqnarray*}
|X_{s}-X_{t_{k-1}}|^2
&\leq &2\left| \int_{t_{k-1}}^{s}(b(X_{u},\theta_{0})-b(X_{t_{k-1}},\theta_{0}))du\right|^2 \\
&\ & +2\left(|b(X_{t_{k-1}},\theta_{0})|(s-t_{k-1})
+\varepsilon |L_{s}-L_{t_{k-1}}|\right)^2 \\
&\leq &2 K^2n^{-1}\int_{t_{k-1}}^{s}|X_{u}-X_{t_{k-1}}|^2du \\
&\ & +2\left(n^{-1}|b(X_{t_{k-1}},\theta_{0})|
+\varepsilon \sup_{t_{k-1}\leq s \leq t_{k}}|L_{s}-L_{t_{k-1}}|\right)^2.
\end{eqnarray*}
By Gronwall's inequality, we get
$$ |X_{s}-X_{t_{k-1}}|^2\leq 2\left(n^{-1}|b(X_{t_{k-1}},\theta_{0})|
+\varepsilon \sup_{t_{k-1}\leq s \leq t_{k}}|L_{s}-L_{t_{k-1}}|\right)^2e^{2K^2n^{-1}(s-t_{k-1})}. $$
It further follows that
\begin{equation}\label{e.3.7}
\sup_{t_{k-1}\leq s \leq t_{k}}|X_{s}-X_{t_{k-1}}|\leq \sqrt{2}\left(n^{-1}|b(X_{t_{k-1}},\theta_{0})|
+\varepsilon \sup_{t_{k-1}\leq s \leq t_{k}}|L_{s}-L_{t_{k-1}}|\right)e^{K^2/n^{2}}
\end{equation}
Thus, by the Lipschitz condition on $b$ and \eref{e.3.7}, we get
\begin{align*}
|H_{n,\varepsilon}^{(1)}(\theta_{0})|
&\leq  \varepsilon^{-1}\sum_{k=1}^{n}|\partial_{\theta_{i}}b(X_{t_{k-1}},\theta_{0})|\cdot
\left|\int_{t_{k-1}}^{t_{k}}(b(X_{s},\theta_{0})-b(X_{t_{k-1}},\theta_{0}))ds\right| \nonumber \\
&\leq  \varepsilon^{-1}\sum_{k=1}^{n}|\partial_{\theta_{i}}b(X_{t_{k-1}},\theta_{0})|\cdot
\int_{t_{k-1}}^{t_{k}}|b(X_{s},\theta_{0})-b(X_{t_{k-1}},\theta_{0})|ds \nonumber \\
&\leq  \varepsilon^{-1}\sum_{k=1}^{n}|\partial_{\theta_{i}}b(X_{t_{k-1}},\theta_{0})|
\int_{t_{k-1}}^{t_{k}}K|X_{s}-X_{t_{k-1}}|ds \nonumber \\
&\leq  (n\varepsilon)^{-1}K\sum_{k=1}^{n}|\partial_{\theta_{i}}b(X_{t_{k-1}},\theta_{0})|\sup_{t_{k-1}\leq s \leq t_{k}}|X_{s}-X_{t_{k-1}}| \nonumber \\
&\leq  \frac{\sqrt{2}Ke^{K^2/n^2}}{n\varepsilon} \cdot \frac{1}{n}\sum_{k=1}^{n}|\partial_{\theta_{i}}b(X_{t_{k-1}},\theta_{0})| \cdot |b(X_{t_{k-1}},\theta_{0})| \nonumber \\
&\quad +\frac{\sqrt{2}Ke^{K^2/n^2}}{n}\sum_{k=1}^{n}|\partial_{\theta_{i}}b(X_{t_{k-1}},\theta_{0})|\sup_{t_{k-1}\leq s \leq t_{k}}|L_{s}-L_{t_{k-1}}| \nonumber \\
&:= H_{n,\varepsilon}^{(1,1)}(\theta_{0})+H_{n,\varepsilon}^{(1,2)}(\theta_{0}).
\end{align*}
It is easy to see that $H_{n,\varepsilon}^{(1,1)}(\theta_{0})$ converges to zero in probability as $n\varepsilon \to
\infty$ since
$$
\frac{1}{n}\sum_{k=1}^{n}|\partial_{\theta_{i}}b(X_{t_{k-1}},\theta_{0})| \cdot |b(X_{t_{k-1}},\theta_{0})|
\leq CK\left(1+\sup_{0\leq s \leq 1}|X_{s}|\right)^{\lambda +1}<\infty \quad  a.s.
$$
{(cf. (\ref{e.2.5})).} By using the basic fact that $$\frac{1}{n}\sum_{k=1}^{n} \sup_{t_{k-1}\leq s \leq t_{k}}|L_{s}-L_{t_{k-1}}|=o_{P}(1),$$
we find that
$$
H_{n,\varepsilon}^{(1,2)}(\theta_{0})\leq \sqrt{2}K e^{K^{2}/n^{2}} C\left(1+\sup_{0\leq s \leq 1}|X_{s}|\right)^{\lambda }
\frac{1}{n}\sum_{k=1}^{n} \sup_{t_{k-1}\leq s \leq t_{k}}|L_{s}-L_{t_{k-1}}|,
$$
which converges to zero in probability as $\varepsilon \to 0$ and $n \to \infty$. Therefore the proof is complete.$\hfill \square $

\begin{lemma}\label{t.3.4}
Assume (A1)-(A4). Then, we have
$$ \sup_{\theta \in \Theta}|K_{n,\varepsilon}(\theta)-K(\theta)|
\overset{P_{\theta_{0}}}{\longrightarrow} 0 $$
as $\varepsilon \to 0$ and $n \to \infty$.
\end{lemma}
{\bf Proof}. It suffices to prove that for $1\leq i, j \leq p$
$$ \sup_{\theta \in \Theta}|K_{n,\varepsilon}^{ij}(\theta)-K^{ij}(\theta)|
\overset{P_{\theta_{0}}}{\longrightarrow} 0 $$
as $\varepsilon \to 0$ and $n \to \infty$.
Note that
\begin{eqnarray*}
K_{n,\varepsilon}^{ij}(\theta)&=&\partial_{\theta_{j}}G_{n,\varepsilon}^{i}(\theta) \nonumber \\
&=&
\sum_{k=1}^{n}(\partial_{\theta_{j}}\partial_{\theta_{i}}
b)^{T}(X_{t_{k-1}},\theta)(X_{t_{k}}-X_{t_{k-1}}-b(X_{t_{k-1}},\theta_{0})\Delta t_{k-1}) \nonumber \\
&\ & +\frac{1}{n}\sum_{i=1}^{n}\left[(\partial_{\theta_{j}}\partial_{\theta_{i}}
b)^{T}(X_{t_{k-1}},\theta)(b(X_{t_{k-1}},\theta_{0})-b(X_{t_{k-1}},\theta)) \right. \nonumber \\
&\ & -\left.(\partial_{\theta_{i}}b)^{T}(X_{t_{k-1}},\theta)\partial_{\theta_{j}}b(X_{t_{k-1}},\theta)\right] \nonumber \\
&:=& K_{n,\varepsilon}^{ij,(1)}(\theta)+K_{n,\varepsilon}^{ij,(2)}(\theta).
\end{eqnarray*}
By using Lemma 3.5 and letting $f(x,\theta)=\partial_{\theta_{j}}\partial_{\theta_{i}}b_{l}(x,\theta)$ ($1\leq i,\ j\leq p, \ 1\leq l \leq d$), we have that
$\sup_{\theta \in \Theta}|K_{n,\varepsilon}^{ij,(1)}(\theta)|$ converges to zero in probability as
$\varepsilon \to 0$ and $n \to \infty$. By using Lemma 3.3 and letting
$f(x,\theta)=(\partial_{\theta_{j}}\partial_{\theta_{i}}b)^{T}(x,\theta)(b(x,\theta_{0})-b(x,\theta))
-(\partial_{\theta_{i}}b)^{T}(x,\theta)\partial_{\theta_{j}}b(x,\theta)$, it follow that
$\sup_{\theta \in \Theta}|K_{n,\varepsilon}^{ij,(2)}(\theta)-K^{ij}(\theta)|$ converges to zero in probability as
$\varepsilon \to 0$ and $n \to \infty$. Thus, the proof is complete. $\hfill \square $ \\

Finally we are ready to prove Theorem \ref{t.3.1}. \\

\smallskip\noindent
{\bf Proof of Theorem \ref{t.3.1}}. The proof ideas mainly follow Uchida \cite{Uch04}.
Let $B(\theta_{0};\rho)=\{\theta: |\theta-\theta_{0}|\leq \rho\}$ for $\rho>0$. Then, by the consistency of $\hat{\theta}_{n,\varepsilon}$, there exists a sequence
$\eta_{n,\varepsilon} \to 0$ as $\varepsilon \to 0$ and $n \to \infty$ such that
$B(\theta_{0};\eta_{n,\varepsilon}) \subset \Theta_0$, and that $P_{\theta_{0}}[\hat{\theta}_{n,\varepsilon} \in B(\theta_{0};\eta_{n,\varepsilon})] \to 1$.
When $\hat{\theta}_{n,\varepsilon} \in B(\theta_{0};\eta_{n,\varepsilon})$, it follows by Taylor's formula that
\[
D_{n,\varepsilon} S_{n,\varepsilon}=\varepsilon^{-1}G_{n,\varepsilon}(\hat{\theta}_{n,\varepsilon})-
\varepsilon^{-1}G_{n,\varepsilon}(\theta_{0}),
\]
where $D_{n,\varepsilon}=\int_{0}^{1} K_{n,\varepsilon}(\theta_{0}+u(\hat{\theta}_{n,\varepsilon}-\theta_{0}))du$
and $S_{n,\varepsilon}=\varepsilon^{-1}(\hat{\theta}_{n,\varepsilon}-\theta_{0})$ {since $B(\theta_{0};\eta_{n,\varepsilon})$ is a convex subset of $\Theta_0$}. We have
\begin{align*}
{|D_{n,\varepsilon}-K_{n,\varepsilon}(\theta_{0})|1_{\{\hat{\theta}_{n,\varepsilon} \in B(\theta_{0};\eta_{n,\varepsilon})\}}}
&\leq  \sup_{\theta \in B(\theta_{0};\eta_{n,\varepsilon})}|K_{n,\varepsilon}(\theta)-K_{n,\varepsilon}(\theta_{0})| \\
& \leq  \sup_{\theta \in B(\theta_{0};\eta_{n,\varepsilon})}|K_{n,\varepsilon}(\theta)-K(\theta)| \\
&\quad +\sup_{\theta \in B(\theta_{0};\eta_{n,\varepsilon})}|K(\theta)-K(\theta_{0})|
+|K_{n,\varepsilon}(\theta_{0})-K(\theta_{0})|.
\end{align*}
Consequently, it follows from Lemma 3.7 that
\[
D_{n,\varepsilon} \overset{P_{\theta_{0}}}{\longrightarrow} K(\theta_{0}),\quad \varepsilon \to 0,\ n \to \infty.
\]
Note that $K(\theta)$ is continuous with respect to $\theta$. Since $-K(\theta_{0})=I(\theta_{0})$ is positive definite, there exists a positive
constant $\delta>0$ such that $\inf_{|w|=1}|K(\theta_{0})w|>2\delta$. For such a $\delta>0$, there exists $\varepsilon(\delta)>0$ and $N(\delta)\in {\mathbb N}$
such that for any $\varepsilon \in (0,\varepsilon(\delta))$, $n>N(\delta)$,
$B(\theta_{0};\eta_{n,\varepsilon}) \subset \Theta_{0}$ and $|K(\theta)-K(\theta_{0})|<\delta/2$ for
$\theta \in B(\theta_{0};\eta_{n,\varepsilon})$,
For such $\delta>0$, let
\[
\Gamma_{n,\varepsilon}=\left\{ \sup_{{|\theta -\theta_0|<\eta_{n,\epsilon}}}|K_{n,\varepsilon}(\theta)-K(\theta_{0})|<\frac{\delta}{2},
\hat{\theta}_{n,\varepsilon} \in B(\theta_{0};\eta_{n,\varepsilon})\right\}.
\]
Then, for any $\varepsilon \in (0,\varepsilon(\delta))$ and $n>N(\delta)$, we have, on $\Gamma_{n,\varepsilon}$,
\begin{eqnarray*}
\sup_{|w|=1}|(D_{n,\varepsilon}-K(\theta_{0}))w|&\leq & \sup_{|w|=1}\left|\left(D_{n,\varepsilon}-\int_{0}^{1}K(\theta_{0}+u(\hat{\theta}_{n,\varepsilon}-\theta_{0}))du\right)w\right|
\\
&\ & +\sup_{|w|=1}\left|\left(\int_{0}^{1}K(\theta_{0}+u(\hat{\theta}_{n,\varepsilon}-\theta_{0}))du-K(\theta_{0})\right)w\right| \\
&\leq & \sup_{|\theta-\theta_{0}|\leq \eta_{n,\varepsilon}}|K_{n,\varepsilon}(\theta)-K(\theta)|+\frac{\delta}{2}<\delta.
\end{eqnarray*}
Thus,
on $ \Gamma_{n,\varepsilon}$,
$$ \inf_{|w|=1}|D_{n,\varepsilon}w|\geq \inf_{|w|=1}|K(\theta_{0})w|-\sup_{|w|=1}|(D_{n,\varepsilon}-K(\theta_{0}))w|>2\delta-\delta=\delta>0. $$
Hence, letting
\[
{{\cal D}_{n,\varepsilon}=\{ D_{n,\varepsilon}\  \text{is invertible, } {\hat{\theta}_{n,\varepsilon} \in B(\theta_{0};\eta_{n,\varepsilon})} \},}
\]
we see that
$P_{\theta_{0}}[{\cal D}_{n,\varepsilon}]\geq P_{\theta_{0}}[\Gamma_{n,\varepsilon}] \to 1$ as $\varepsilon \to 0$ and $n \to \infty$ by Lemma 3.7.
Now set
\[
{ U_{n,\varepsilon}= D_{n,\varepsilon} 1_{{\cal D}_{n,\varepsilon}} + I_{p\times p} 1_{{\cal D}_{n,\varepsilon}^c}},
\]
where $I_{p\times p}$ is the identity matrix. Then it is easy to see that
$$ |U_{n,\varepsilon}-K(\theta_{0})|\leq |D_{n,\varepsilon}-K(\theta_{0})|1_{{\cal D}_{n,\varepsilon}}+
|I_{p\times p}-K(\theta_{0})| 1_{{\cal D}_{n,\varepsilon}^{c}}\overset{P_{\theta_{0}}}{\longrightarrow} 0, $$
since $P_{\theta_{0}}[{\cal D}_{n,\varepsilon}] \to 1$.
Thus, by Lemma 3.6, we obtain that
\begin{align*}
S_{n,\varepsilon}
&=  U_{n,\varepsilon}^{-1}D_{n,\varepsilon}S_{n,\varepsilon}1_{{\cal D}_{n,\varepsilon}} +  S_{n,\varepsilon}1_{{\cal D}_{n,\varepsilon}^c}\\
&= U_{n,\varepsilon}^{-1}(-\varepsilon^{-1}G_{n,\varepsilon}(\theta_{0}))1_{{\cal D}_{n,\varepsilon}} +  S_{n,\varepsilon}1_{{\cal D}_{n,\varepsilon}^c}\\
&\overset{P_{\theta_{0}}}{\longrightarrow}  (I(\theta_{0}))^{-1}\left(\int_{0}^{1}(\partial_{\theta_{1}}b)^{T}(X_{s}^{0},\theta_{0})dL_{s},
\dots, \int_{0}^{1}(\partial_{\theta_{p}}b)^{T}(X_{s}^{0},\theta_{0})dL_{s}\right)^{T}
\end{align*}
as  $\varepsilon \to 0$, $n \to \infty$ and $n \varepsilon \to \infty$.
This completes the proof. $\hfill \square  $ \\

\end{document}